\documentclass[fleqn,11pt]{article}

\usepackage{amsmath, amscd, amsfonts, amssymb, graphicx, color}
\usepackage{hyperref,dashrule}
\usepackage{cite}
\usepackage{amsmath}
\usepackage{amsfonts}
\usepackage{epsfig}
\usepackage{amssymb}
\usepackage{amsthm}
\usepackage{epstopdf}
\usepackage{footnote}
\usepackage{newlfont}
\usepackage{color}
\newtheorem{theorem}{Theorem}[section]
\newtheorem{lemma}[theorem]{Lemma}

\newtheorem{corollary}[theorem]{Corollary}
\theoremstyle{definition}
\newtheorem{definition}[theorem]{Definition}

\newtheorem{remark}[theorem]{Remark}

\begin{document}

\title{Prescribed Eigenvalues via Optimal Perturbation of main-diagonal submatrix}

\author{M. R. Eslahchi\thanks{Department of Computer Science, Faculty of Mathematical Sciences, Tarbiat Modares University, P.O. Box 14115-134, Tehran, Iran (eslahchi@modares.ac.ir).}\,
and E. Kokabifar \thanks{Department of Mathematics Education, Farhangian University, P.O. Box 14665-889, Tehran, Iran (kokabifar@cfu.ac.ir).\newline These authors were supported by grant no. 4033159
of Iran National Science Foundation (INSF).\newline Corresponding author: M. R. Eslahchi.}}
\date{}
\maketitle

\vspace{-6mm}

\begin{abstract}
Consider a given square matrix $\textrm {K}$ with square blocks $A_{11},A_{22},\ldots,A_{nn}$ on the main diagonal. This paper aims to compute an optimal perturbation $\Delta$ of a preassigned block $A_{ii}\in\mathbb{C}^{d_i\times d_k}, \left(1\le i\le n\right)$,with respect to the spectral norm distance, such that the perturbed matrix ${\textrm {K}_X}$ has $k \le d_i$ prescribed eigenvalues.
This paper presents a method for constructing  the optimal perturbation by improving and extending the methodology, necessary definitions and lemmas of previous related works. Some conceivable applications of this subject are also presented. Numerical
experiments are provided to illustrate the validity of the method.

\end{abstract}



\section{Introduction}\label{intr}

Consider $A\in\mathbb{C}^{n\times n}$ with distinct eigenvalues. Finding the nearest matrix $X$ to $A$ with multiple eigenvalues is known as Wilkinson's problem \citeonline{wilkinson}. There are quite a few works dealing with the problem, see for instance \citeonline{wil1,wil4,demmel1} in which some bounds for the distance are computed. Anyway, as far as the authors are aware there is no explicit solution for Wilkinson's eigenvalue problem. It is worth noting that Malyshev \cite{malyshev} computed the spectral norm distance from $A$ to all $n\times n$ matrices having a \textit{prescribed multiple} eigenvalue. Even though it is often observed that, for a specific eigenvalue, matrices with algebraic multiplicity 2 are closer to $A$ rather than those having greater algebraic multiplicities \citeonline{wilkinson}, results obtained in \cite{malyshev} are still under a strict condition: A \textit{predefined} multiple eigenvalue.

Lippert \cite{lipert} and Gracia \cite{gracia} constructed the optimal perturbation of $A$ such that the perturbed matrix has two prescribed eigenvalues. Some other related works are \cite{lipertk,mengi,ikramov}. Specifically, Kokabifar et al. \cite{klk}, motivated by the above spectrum updating problems, computed the spectral norm distance from  a given matrix $A\in\mathbb{C}^{n\times n}$ to the set of $n\times n$ matrices that have $k \le n$ prescribed distinct eigenvalues. What if we  are allowed to change only one part of a given matrix in order to make prescribed eigenvalues? To answer this question, Gracia and Velasco\cite{graciavelasco}, partitioned a given square matrix $\textrm {K}$ as 
\begin{equation}\label{matrixk}
\textrm {K}=\left[ {\begin{array}{*{20}{c}}
A_{n\times n}&B_{n\times m} \\ 
C_{m\times n}&D_{m\times m} 
\end{array}} \right], \qquad \textnormal{for some } m,n\in\mathbb{N},
\end{equation}
made optimal changes to $D$, then making good use of Theorem \ref{mainrank3} found matrix $\Delta$ with minimum 2-norm such that the perturbed matrix $\textrm {K}_X=\left[ {\begin{array}{*{20}{c}}
A&B \\ 
C&D+\Delta 
\end{array}} \right]=\left[ {\begin{array}{*{20}{c}}
A&B \\ 
C&X 
\end{array}} \right]$ has a multiple a prescribed eigenvalue. This problem is also addressed for two given eigenvalues in \cite{tokhmesag}. There are some other papers in which only the same square southeast submatrix of $\textrm {K}$ is perturbed. See \cite{armen,gonz,graciavelasco2} and references therein.

Suppose the matrix K$\in\mathbb{C}^{d\times d}$ such that $A_{11},A_{22},\ldots,A_{nn}, \left(A_{ii}\in\mathbb{C}^{d_i\times d_i}\right)$ are its main diagonal blocks and a set of given complex numbers $\Lambda  = \{ \lambda _1 ,\lambda _2 , \ldots ,\lambda _k \}$. Denote by $\mathcal{M}_k$, $d$-square matrices that $\Lambda$ is a subset of their spectrum.

This paper describes a clear and intelligible computational technique to construct the perturbation $\Delta$ such that for preassigned main diagonal block $A_{ii}$, with $k\le d_i$, ${\left\| {\Delta} \right\|_2}$  has the minimum spectral norm and K$_X$ belongs to $\mathcal{M}_k$, when intersection of spectrum of $A_{jj}, (j\neq i)$ and $\Lambda$ is empty. Expanding and improving the methodology used in \cite{klk,graciavelasco,lipertk,lipert,tokhmesag}, presenting a general solution for the matrix nearness problem in an evident computational manner as well as introducing a number of possible applications of the problem are the main goals considered herein.

In the first place, for the sake of simplicity, we consider the problem when only three prescribed eigenvalues are given and K is partitioned as in (\ref{matrixk}), i.e., only the southeast submatrix is perturbed to make predefined eigenvalues. Calculations and proofs presented in this section will clarify those of next sections in which the ideas and computations are generalized for $k\le d_n$ given scalars. Finally, these results are extended to every main diagonal block, benefiting from permutation matrices.

\section{Nearest submatrix providing three prescribed eigenvalues}\label{sec2}
In this section, we start with calculating a lower bound for 2-norm of the perturbation $\Delta$. See Lemma \ref{lemmalowerbound3}. Then, in Subsection \ref{cons3}, an optimal perturbation with minimum possible norm is constructed.
\subsection{Computing lower bounds for the optimal perturbation}\label{lowe}

For a given set $\Lambda  = \{ \lambda _1 ,\lambda _2 , \lambda _3 \}$ and a square matrix $T\in\mathbb{C}^{l\times l}$, define 
\[{Q_T}(\gamma_{12},\gamma_{13},\gamma_{23} ) = \left[ {\begin{array}{*{20}{c}}
{T - {\lambda _1}{I_l}}&{{\gamma _{12}}{I_l}}&{{\gamma _{13}}{I_l}} \\ 
0&{T - {\lambda _2}{I_l}}&{{\gamma _{23}}{I_l}} \\ 
0&0&{T - {\lambda _3}{I_l}} 
\end{array}} \right],\qquad \gamma_{12},\gamma_{13},\gamma_{23} \in \mathbb{C}, \]
and denote $\gamma_{12},\gamma_{13},\gamma_{23} \in \mathbb{C}$ by $\gamma$. Similar to the proof of Lemma 2.1 of \cite{klk} (see also Lemma 7 of \cite{lipert}) it can be concluded that if $T$ has $\lambda _1 ,\lambda _2 ,\lambda _3$ as some of its eigenvalues, then for all $\gamma$ it holds that ${s_{3l-2} }\left( {Q_T(\gamma )} \right)=0.$ More importantly, the following corollary can be derived.
\begin{corollary}\label{colorank3}
Let $\lambda _1 ,\lambda _2 , \lambda _3$ be eigenvalues of $T \in {\mathbb{C}^{d \times d}}$. Then for all $\gamma$  we have
\begin{equation}\label{coro}
{\textnormal{rank}\left( {\left[ {\begin{array}{*{20}{c}}
{T - {\lambda _1}{I_l}}&{{\gamma _{12}}{I_l}}&{{\gamma _{13}}{I_l}} \\ 
0&{T - {\lambda _2}{I_l}}&{{\gamma _{23}}{I_l}} \\ 
0&0&{T - {\lambda _3}{I_l}} 
\end{array}} \right]} \right) \leqslant 3d - 3.}
\end{equation}
\end{corollary}
The next theorem, which is indeed Theorem 1.1 of \cite{gonz} as well as Theorem 5 of \cite{graciavelasco} is of great importance for the rest of discussions. 

\begin{theorem}\label{mainrank3}
Recall the matrix $\textnormal{K}$ as in (\ref{matrixk}) and let $ \textnormal {K}_X=\left[ {\begin{array}{*{20}{c}}
A_{n\times n}&B_{n\times m} \\ 
C_{m\times n}&X_{m\times m} 
\end{array}} \right]$,
\[\begin{gathered}
\rho  = \textnormal{rank}\left( {\left[ {A,B} \right]} \right) + \textnormal{rank}\left( {\left[ \begin{gathered}
A \hfill \\
B \hfill \\ 
\end{gathered}  \right]} \right) - \textnormal{rank}\left( A \right), \hfill \\
\begin{array}{*{20}{c}}
{M = \left( {I - A{A^\dag }} \right)B,}&{N = C\left( {I - {A^\dag }A} \right)},&{S\left( X  \right) = \left( {I - N{N^\dag }} \right)\left( {X  - C{A^\dag }B} \right)\left( {I - {M^\dag }M} \right)},
\end{array} \hfill \\ 
\end{gathered} \]
where $A^\dag$ denotes the\textit{ Moore–Penrose} pseudo-inverse of $A$, then for every $X\in \mathbb{C}^{m\times m}$, it holds that
\[\textnormal{rank}\left( {{\textnormal {K}_X }} \right) = \rho  + \textnormal{rank}\left( {S\left( X  \right)} \right),\]
furthermore, If $r$ is an integer satisfying the  inequalities $\rho  \leqslant r < \textnormal {rank}\left(\textnormal {K} \right)$, then
\[\min \left\{ {\left\| {X  - D} \right\|:X  \in {\mathbb{C}^{m \times m}},\textnormal{rank}\left( {{\textnormal {K}_X }} \right) \leqslant r} \right\} = {s_{r - \rho  + 1}}\left( {S\left( X  \right)} \right).\]
\end{theorem}
To investigate the properties of the singular values of a matrix-valued function and its associated singular vectors, results deduced in Lemma 7 of \cite{graciavelasco} (see also Theorem 3.4 of \cite{klk} for more details) are required.
\begin{theorem}\label{analytic}
Let $\Omega$ be an open subset of $\mathbb{R}$, and let $F:\Omega  \to {\mathbb{C}^{m \times n}}$ be an analytic function on $\Omega$. If the function ${s_i}\left( {F\left( t \right)} \right)$ has a positive local extremum at $t_0 \in \Omega$, then there exists a pair of singular vectors $u\in \mathbb{C}^{m\times 1},v\in \mathbb{C}^{n\times 1}$ of $F\left(t_0\right)$ corresponding to ${s_i}\left( {F\left( t_0\right)} \right)$ such that $\operatorname{Re} \left( {{u^*}\dfrac{{dF}}{{dt}}\left( {{t_0}} \right)v} \right) = 0.$
\end{theorem}
Let
\[\begin{array}{*{20}{c}}
  {\mathcal{A} = \left[ {\begin{array}{*{20}{c}}
  {A - {\lambda _1}{I_n}}&{{\gamma _{12}}{I_n}}&{{\gamma _{13}}{I_n}} \\ 
  0&{A - {\lambda _2}{I_n}}&{{\gamma _{23}}{I_n}} \\ 
  0&0&{A - {\lambda _3}{I_n}} 
\end{array}} \right],}&{\mathcal{B} = \left[ {\begin{array}{*{20}{c}}
  B&0&0 \\ 
  0&B&0 \\ 
  0&0&B 
\end{array}} \right],} \\ 
  {\mathcal{C} = \left[ {\begin{array}{*{20}{c}}
  C&0&0 \\ 
  0&C&0 \\ 
  0&0&C 
\end{array}} \right],}&{\mathcal{X} = \left[ {\begin{array}{*{20}{c}}
  {X - {\lambda _1}{I_m}}&{{\gamma _{12}}{I_m}}&{{\gamma _{13}}{I_m}} \\ 
  0&{X - {\lambda _2}{I_m}}&{{\gamma _{23}}{I_m}} \\ 
  0&0&{X - {\lambda _3}{I_m}} 
\end{array}} \right],} 
\end{array}\]
if $\lambda_1,\lambda_2$ and $\lambda_3$ belong to spectrum of $ \textnormal {K}_X$,
 so Corollary \ref{colorank3} concludes
\begin{eqnarray*}
\textnormal {rank}\left( \left[{\begin{array}{*{20}{c}}
\mathcal A&\mathcal B \\ 
\mathcal	C&\mathcal X 
\end{array}}\right] \right) &=&\textnormal {rank}\left( {\left[ {\begin{array}{*{20}{c}}
{{K_X}  - {\lambda _1}{I_{m + n}}}&{{\gamma _{12}}{I_{m + n}}}&{{\gamma _{13}}{I_{m + n}}} \\ 
0&{{K_X} - {\lambda _2}{I_{m + n}}}&{{\gamma _{23}}{I_{m + n}}} \\ 
0&0&{{K_X} - {\lambda _m}{I_{m + n}}} 
\end{array}} \right]} \right) \\&\leqslant& 3\left( {n + m} \right)-3,
\end{eqnarray*}
then, by applying Theorem \ref{mainrank3}, we have
$\textnormal {rank}\left( \left[{\begin{array}{*{20}{c}}
\mathcal A&\mathcal B \\ 
\mathcal	C&\mathcal X 
\end{array}}\right] \right) =  \rho \left( \gamma  \right) + \textnormal {rank}\left( {\mathcal{S}_3\left( X ,\gamma\right)} \right),$
in which
\begin{equation}\label{mains3}
\begin{gathered}
\mathcal{M}\left( \gamma  \right) = \left( {{I_{3n}} - \mathcal{A}{\mathcal{A}^\dag }} \right)\mathcal{B},\qquad \mathcal{N}\left( \gamma  \right) = \mathcal{C}\left( {{I_{3n}} - {\mathcal{A}^\dag }\mathcal{A}} \right), \hfill \\
\mathcal{S}_3\left( X ,\gamma\right) = \left( {I - \mathcal{N}{\mathcal{N}^\dag }} \right)\left( {\mathcal{X} - \mathcal{C}{\mathcal{A}^\dag }\mathcal{B}} \right)\left( {I - {\mathcal{M}^\dag }\mathcal{M}} \right), \hfill \\ 
\end{gathered}
\end{equation}
and
\begin{equation}\label{rho}
\rho\left(\gamma\right)  = \textnormal{rank}\left( {\left[ {\mathcal{A},\mathcal{B}} \right]} \right) + \textnormal{rank}\left( {\left[ \begin{gathered}
\mathcal{A} \hfill \\
\mathcal{B} \hfill \\ 
\end{gathered}  \right]} \right) - \textnormal{rank}\left( \mathcal{A} \right)=3n,
\end{equation}	
thus, $\rho \left( \gamma  \right) + \textnormal{rank}\left( {{\mathcal{S}_3}\left( {X,\gamma } \right)} \right) \leqslant 3\left( {n + m - 1} \right)$ implies $\textnormal{rank}\left( {{\mathcal{S}_3}\left( {X,\gamma } \right)} \right) \leqslant 3\left(m-1\right)$
and consequently for every $\gamma_{12},\gamma_{13},\gamma_{23} \in \mathbb{C}$, it holds that
${s_{3m-2}}\left( {{\mathcal{S}_3}\left( {X,\gamma } \right)} \right) = 0.$

Weyl inequalities for singular values (e.g. see Corollary5.1 of\cite{demmel})  and the last relation provide us the lower bound we were looking for. The proof of the following lemma is analogous to Lemma 22 of \cite{graciavelasco} (see also Lemma 2.3 of \cite{klk}).

\begin{lemma}\label{lemmalowerbound3}
Assume that the matrix $\textnormal{K}$ as in (\ref{matrixk}) and the set $\Lambda  = \{ \lambda _1 ,\lambda _2 , \lambda _3 \}$ are given. If for $X\in \mathbb{C}^{m\times m}$, spectrum of $\textnormal {K}_X=\left[ {\begin{array}{*{20}{c}}
A&B \\ 
C&X 
\end{array}} \right]$ includes $\Lambda$ then
\[{\left\| {X - D} \right\|_2}\ge {s_{3m-2}}\left( {{\mathcal{S}_3}\left( {D,\gamma } \right)} \right).\]
\end{lemma}

\subsection{Construction of the optimal perturbation}\label{cons3}
This section concerns computing an optimal perturbation $\Delta_*$ satisfying ${\left\| \Delta_* \right\|_2}= s_{3m-2} \left( \gamma  \right)$ for some value of $\gamma$, and K$_{X_*}\in \mathcal{M}_3$. To do this, let $A^{-1}_i$ denotes $\left( {A - {\lambda _i}{I_n}} \right)^{ - 1}$, then
\[{\mathcal{A}^\dag } = \left[ {\begin{array}{*{20}{c}}
{A_1^{ - 1}}&{ - {\gamma _{12}}A_1^{ - 1}A_2^{ - 1}}&{ - {\gamma _{13}}A_1^{ - 1}A_3^{ - 1} + {\gamma _{12}}{\gamma _{23}}A_1^{ - 1}A_2^{ - 1}A_3^{ - 1}} \\ 
0&{A_2^{ - 1}}&{ - {\gamma _{23}}A_2^{ - 1}A_3^{ - 1}} \\ 
0&0&{A_3^{ - 1}} 
\end{array}} \right],\]
and so, using (\ref{mains3}) we have
\begin{eqnarray*}
&&\hspace{-.7cm}{\mathcal{S}_3}\left( {D,\gamma } \right) =\\&&\hspace{-.7cm} \small{\left[ {\begin{array}{*{20}{c}}
{{D_1} - CA_1^{ - 1}B}&{{\gamma _{12}}\left( {{I_m} + CA_1^{ - 1}A_2^{ - 1}B} \right)}&{{\gamma _{13}}\left( {{I_m} + CA_1^{ - 1}A_3^{ - 1}B} \right) - {\gamma _{12}}{\gamma _{23}}CA_1^{ - 1}A_2^{ - 1}A_3^{ - 1}B} \\ 
0&{D_2 - CA_2^{ - 1}B}&{{\gamma _{23}}\left( {{I_m} + CA_2^{ - 1}A_3^{ - 1}B} \right)} \\ 
0&0&{D_3 - CA_3^{ - 1}B} 
\end{array}} \right]},
\end{eqnarray*}
where $D_i =\left( {D - {\lambda _i}{I_n}} \right)$. Introducing
\begin{equation}\label{MNP}
\begin{array}{*{20}{c}}
{{M_i} = {D_i} - CA_i^{ - 1}B,}&{{N_{ij}} = {I_m} + CA_i^{ - 1}A_j^{ - 1}B,}&{P_{123} = CA_1^{ - 1}A_2^{ - 1}A_3^{ - 1}B,} 
\end{array}
\end{equation}
$\mathcal{S}_3\left( {D,\gamma } \right)$ reduces to
\[\mathcal{S}_3\left( {D,\gamma } \right)=\left[ {\begin{array}{*{20}{c}}
{{M_1}}&{{\gamma _{12}}{N_{12}}}&{{\gamma _{13}}{N_{13}} - {\gamma _{12}}{\gamma _{23}}P_{123}} \\ 
0&{{M_2}}&{{\gamma _{23}}{N_{23}}} \\ 
0&0&{{M_3}} 
\end{array}} \right].\]

In what follows, we obtain further properties of ${s_{3m-2} }\left( {Q_A(\gamma )} \right)$ and its associated singular vectors. In the next section, this properties are applied to construct the optimal perturbation $\Delta_*$. For our discussion, it is necessary to reform some definitions and lemmas of \cite{malyshev, gracia, lipert,mengi}.
Suppose that vectors
\begin{equation*}
u(\gamma)=\left[ {\begin{array}{*{20}{c}}
{{u_1}(\gamma )}\\{{u_2}(\gamma )}\\
{{u_3}(\gamma )}
\end{array}} \right], v(\gamma)=\left[ {\begin{array}{*{20}{c}}
{{v_1}(\gamma )}\\{{v_2}(\gamma )}\\
{{v_3}(\gamma )}
\end{array}} \right] \in {\mathbb{C}^{3m\times1}}~({u_j}(\gamma ),{v_j}(\gamma ) \in {\mathbb{C}^{m\times1}},j = 1,2,3),
\end{equation*}
is a pair of left and right singular vectors of ${s_{3m-2} }\left( \gamma \right)$, respectively. Define
\begin{equation*}
U(\gamma) = [{u_1}(\gamma ), {u_2}(\gamma ) ,{u_3}(\gamma )]_{m \times 3},\qquad {  \mbox{and}}\qquad V(\gamma) = [{v_1}(\gamma ),{v_2}(\gamma ) ,{v_3}(\gamma )]_{m \times 3}.
\end{equation*}

Considering definition of the vectors $u(\gamma)$ and $v(\gamma)$, we have 
\begin{equation*}\label{a}
{\mathcal{S}_3}\left( {D,\gamma } \right)v(\gamma ) = {s_{3m-2} }\left(\gamma \right)u(\gamma ),
\end{equation*}
or equivalently 
\begin{equation}\label{setofeq}
\left\{ {\begin{array}{*{20}{c}}
  {{M_1}{v_1}\left( \gamma  \right) + {\gamma _{12}}{N_{12}}{v_2}\left( \gamma  \right) + \left( {{\gamma _{13}}{N_{13}} - {\gamma _{12}}{\gamma _{23}}P_{123}} \right){v_3}\left( \gamma  \right) = s{u_1}\left( \gamma  \right)}, \\ 
  {{M_2}{v_2}\left( \gamma  \right) + {\gamma _{23}}{N_{23}}{v_3}\left( \gamma  \right) = s{u_2}\left( \gamma  \right)}, \\ 
  {{M_3}{v_3}\left( \gamma  \right) = s{u_3}\left( \gamma  \right)}. 
\end{array}} \right.
\end{equation}
also without loss of generality, assume that $u(\gamma)$ and $v(\gamma)$ are unit vectors. The inequality for singular values of the sum of two matrices reported in Theorem 3.3.16 of \cite{horn} yields
\begin{eqnarray*}
{s_{3m - 2}}\left( \mathcal{S}_3\left( {D,\gamma } \right) \right) &\leqslant& {s_1}\left( {\left[ {\begin{array}{*{20}{c}}
{{M_1}}&0&0 \\ 
0&{{M_2}}&0 \\ 
0&0&{{M_3}} 
\end{array}} \right]} \right) \\ &+& {s_{3m - 2}}\left( {\left[ {\begin{array}{*{20}{c}}
0&{{\gamma _{12}}{N_{12}}}&{{\gamma _{13}}{N_{13}} - {\gamma _{12}}{\gamma _{23}}P_{123}} \\ 
0&0&{{\gamma _{23}}{N_{23}}} \\ 
0&0&0 
\end{array}} \right]} \right) \\&=& {s_1}\left( {\left[ {\begin{array}{*{20}{c}}
{{M_1}}&0&0 \\ 
0&{{M_2}}&0 \\ 
0&0&{{M_3}} 
\end{array}} \right]} \right) \leqslant \max \left\{ {\left\| {{M_1}} \right\|,\left\| {{M_2}} \right\|,\left\| {{M_3}} \right\|} \right\},
\end{eqnarray*}
furthermore,
\begin{eqnarray*}
{s_{3m - 2}}\left( \mathcal{S}_3\left( {D,\gamma } \right) \right) &=& \frac{1}{{{s_3}\left( {{{\left[ {\begin{array}{*{20}{c}}
{{M_1}}&{{\gamma _{12}}{N_{12}}}&{{\gamma _{13}}{N_{13}} - {\gamma _{12}}{\gamma _{23}}P_{123}} \\ 
0&{{M_2}}&{{\gamma _{23}}{N_{23}}} \\ 
0&0&{{M_3}} 
\end{array}} \right]}^{ - 1}}} \right)}} \\&=& \frac{1}{{{s_3}\left( {\left[ {\begin{array}{*{20}{c}}
{M_1^{ - 1}}&{ - {\gamma _{12}}M_1^{ - 1}{N_{12}}M_2^{ - 1}}&W \\ 
0&{M_2^{ - 1}}&{ - {\gamma _{23}}M_2^{ - 1}{N_{23}}M_3^{ - 1}} \\ 
0&0&{M_3^{ - 1}} 
\end{array}} \right]} \right)}}
\end{eqnarray*}
in which, $W = M_1^{ - 1}\left( { - {\gamma _{13}}{N_{13}}M_3^{ - 1} + {\gamma _{12}}{\gamma _{23}}P_{123}M_3^{ - 1} + {\gamma _{12}}{\gamma _{23}}{N_{12}}M_2^{ - 1}{N_{23}}M_3^{ - 1}} \right)$. Again, using Theorem 3.3.16 of \cite{horn}, we have
\begin{eqnarray*}
&&{s_3}\left( {\left[ {\begin{array}{*{20}{c}}
  {M_1^{ - 1}}&{ - {\gamma _{12}}M_1^{ - 1}{N_{12}}M_2^{ - 1}}&W \\ 
  0&{M_2^{ - 1}}&{ - {\gamma _{23}}M_2^{ - 1}{N_{23}}M_3^{ - 1}} \\ 
  0&0&{M_3^{ - 1}} 
\end{array}} \right]} \right)\\ &\geqslant& {s_3}\left( {\left[ {\begin{array}{*{20}{c}}
  0&{ - {\gamma _{12}}M_1^{ - 1}{N_{12}}M_2^{ - 1}}&W \\ 
  0&0&{ - {\gamma _{23}}M_2^{ - 1}{N_{23}}M_3^{ - 1}} \\ 
  0&0&0 
\end{array}} \right]} \right) - {s_1}\left( {\left[ {\begin{array}{*{20}{c}}
  {M_1^{ - 1}}&0&0 \\ 
  0&{M_2^{ - 1}}&0 \\ 
  0&0&{M_3^{ - 1}} 
\end{array}} \right]} \right),
\end{eqnarray*}
so, if $\left| \gamma  \right| \to \infty$, then ${s_{3}}\left( \mathcal{S}_3\left( {D,\gamma } \right)^{-1} \right) \to \infty$ and consequently ${s_{3m - 2}}\left( \mathcal{S}_3\left( {D,\gamma } \right) \right)\to0$.

Above discussions conclude that there exists a finite point $\gamma $ where the bounded function ${s_{3m-2} }\left( \gamma \right)$ attains its maximum value.

\begin{definition}\label{sstar3}
Let $\gamma_*>0$ be a point where the singular value ${s_{3m-2} }\left(\gamma \right)$ attains its maximum value. We set $\alpha_3^* ={s_{3m-2} }\left( \gamma_* \right).$
\end{definition}
If $\alpha_3^* =0$, then straightforward calculations conclude that $\lambda_1, \lambda_2,\lambda_3$ are some eigenvalues of ${\mathcal{S}_3}\left( {D,\gamma_* } \right)$. Accordingly, hereinafter it is assumed that $\alpha_3^*>0$. In addition, consider $\alpha_3^*$ as a simple singular value of ${\mathcal{S}_3}\left( {D,\gamma_* } \right)$. So, using Theorem \ref{analytic} and undemanding computations analogous to Lemma 3.5 and Lemma 3.7 of \cite{klk} and Lemma 2.7 of \cite{tokhmesag}, result in the next lemma which plays a key role to show  optimality of the perturbation.
\begin{lemma}\label{lem3}
Let $\gamma_*$ and $\alpha_3^*$ be as defined in Definition \ref{sstar3}, $\alpha_3^*>0$ be a simple singular value of ${\mathcal{S}_3}\left( {D,\gamma_*} \right)$, and two matrices $U\left( {{\gamma _*}} \right)$ and $V\left( {{\gamma _*}} \right)$ have full rank. Then
\[U{\left( {{\gamma _*}} \right)^*}U\left( {{\gamma _*}} \right) = V{\left( {{\gamma _*}} \right)^*}V\left( {{\gamma _*}} \right).\]
\end{lemma}

Now, define the desired perturbation as follows
\begin{equation}\label{delta}
\Delta_*  =  - {\alpha_3^*}U({\gamma _*}){V }({\gamma _*})^\dag,
\end{equation}
which implies that
\begin{equation}\label{deltav}
\Delta_* V({\gamma _*}) =  - \alpha _3^*U({\gamma _*}) \Leftrightarrow \Delta_* {v_i}({\gamma _*}) =  - \alpha _3^*{u_i}({\gamma _*}),\qquad i = 1,2,3.
\end{equation}
Lemma \ref{lem3} concludes that  $U(\gamma _* )$ and $V(\gamma _* )$ have the same nonzero singular values. Consequently, a  unitary matrix $W$ exits such that $U(\gamma_*)=W V(\gamma_*)$. Moreover, From Lemma \ref{lem3} one can deduced $ V{(\gamma )^\dag } V{(\gamma )}=I_3$. Thus, 
\begin{equation*}
{\left\| \Delta_*  \right\|_2} = {\left\| { - {\alpha_3^*}U({\gamma _*}){V }({\gamma _*})^\dag} \right\|_2} =\alpha_3^*{\left\| {WV({\gamma _*}){V }({\gamma _*})^\dag} \right\|_2} = \alpha_3^*,
\end{equation*}
which means that the perturbation $\Delta_*$ introduced in (\ref{delta}) meets the optimality condition as stated in Lemma \ref{lemmalowerbound3}. 

Eventually, it is shown that prespecified scalars $\lambda_1,\lambda_2$ and $\lambda_3$ belong to spectrum of $\textnormal {K}_{X_*}=\left[ {\begin{array}{*{20}{c}}
A&B \\ 
C&X_* 
\end{array}} \right]$, in which $X_*=D+\Delta_*$. To this end, keeping in mind rank$\left(V(\gamma_*)\right)=3$, vectors $w_i(\gamma)\in \mathbb{C}^{n\times1}, (i=1,2,3)$ are introduced satisfying
\begin{equation}\label{meq}
\textnormal {K}_{X_*}\left[ {\begin{array}{*{20}{c}}
  {{w_3}\left( \gamma  \right)}&{{w_2}\left( \gamma  \right)}&{{w_1}\left( \gamma  \right)} \\ 
  {{v_3}\left( \gamma  \right)}&{{v_2}\left( \gamma  \right)}&{{v_1}\left( \gamma  \right)} 
\end{array}} \right] = \left[ {\begin{array}{*{20}{c}}
  {{w_3}\left( \gamma  \right)}&{{w_2}\left( \gamma  \right)}&{{w_1}\left( \gamma  \right)} \\ 
  {{v_3}\left( \gamma  \right)}&{{v_2}\left( \gamma  \right)}&{{v_1}\left( \gamma  \right)} 
\end{array}} \right]\left[ {\begin{array}{*{20}{c}}
  {{\lambda _1}}&{ - {\gamma _{23}}}&{ - {\gamma _{13}}} \\ 
  0&{{\lambda _2}}&{ - {\gamma _{12}}} \\ 
  0&0&{{\lambda _3}} 
\end{array}} \right].
\end{equation}
Multiplying the first row of the above matrix equation yields
\[\begin{array}{*{20}{c}}
  {{w_3}\left( \gamma  \right) =  - A_3^{ - 1}B{v_3}\left( \gamma  \right),} \\ 
  {{w_2}\left( \gamma  \right) = A_2^{ - 1}\left( { - B{v_2}\left( \gamma  \right) + {\gamma _{23}}A_3^{ - 1}B{v_3}\left( \gamma  \right)} \right),} \\ 
  {{w_1}\left( \gamma  \right) = A_1^{ - 1}\left( { - B{v_1}\left( \gamma  \right) + {\gamma _{12}}A_2^{ - 1}B{v_2}\left( \gamma  \right) + {\gamma _{13}}A_3^{ - 1}B{v_3}\left( \gamma  \right) - {\gamma _{12}}{\gamma _{23}}A_2^{ - 1}A_3^{ - 1}B{v_3}\left( \gamma  \right)} \right),} 
\end{array}\]
moreover, the three vectors $w_i(\gamma), (i=1,2,3)$ computed from the first row, are also satisfy the second row of the matrix equation (\ref{meq}). To be more precise, multiplying the second row by the first, column obtains
\[C{w_3}\left( \gamma  \right) + {X_*}{v_3}\left( \gamma  \right) = {\lambda _3}{v_3}\left( \gamma  \right),\]
by substituting $w_3(\gamma)$ into the above equation and considering the relation (\ref{deltav}), we have
\[  C\left( { - A_3^{ - 1}B{v_3}\left( \gamma  \right)} \right) + D{v_3}\left( \gamma  \right) - \alpha _3^*{u_3}\left( \gamma  \right) = {\lambda _3}{v_3}\left( \gamma  \right),\]
or 
\[
  \left( {\left( {D - {\lambda _3}{I_m}} \right) - CA_3^{ - 1}B} \right){v_3}\left( \gamma  \right) = \alpha _3^*{u_3}\left( \gamma  \right), \]
which is true according to the last equation of (\ref{setofeq}). Additionally, for the second column we have the following relation
\[C{w_2}\left( \gamma  \right) + {X_*}{v_2}\left( \gamma  \right) =  - {\gamma _{23}}{v_2}\left( \gamma  \right) + {\lambda _2}{v_2}\left( \gamma  \right),\]
that can be written as
\[C\left( { - A_2^{ - 1}B{v_2}\left( \gamma  \right) + {\gamma _{23}}A_2^{ - 1}A_3^{ - 1}B{v_3}\left( \gamma  \right)} \right) + D{v_2}\left( \gamma  \right) - \alpha _3^*{u_2}\left( \gamma  \right) =  - {\gamma _{23}}{v_2}\left( \gamma  \right) + {\lambda _2}{v_2}\left( \gamma  \right),\]
or equivalently
\[\left( {\left( {D - {\lambda _2}{I_m}} \right) - CA_2^{ - 1}B} \right){v_2}\left( \gamma  \right) + {\gamma _{23}}\left( {{I_m} + CA_2^{ - 1}A_3^{ - 1}B} \right){v_3}\left( \gamma  \right) = \alpha _3^*{u_2}\left( \gamma  \right)\]
that is indeed, the second equation of (\ref{setofeq}). Finally, multiplying the second row of the matrix equation (\ref{meq}) by the third column, one can obtain
\[C{w_1}\left( \gamma  \right) + {X_*}{v_1}\left( \gamma  \right) =  - {\gamma _{13}}{v_1}\left( \gamma  \right) - {\gamma _{12}}{v_2}\left( \gamma  \right) + {\lambda _1}{v_1}\left( \gamma  \right),\]
substituting $w_1(\gamma)$ and rearranging the coefficients result in 
\begin{eqnarray*}
&&\left( {\left( {D - {\lambda _1}{I_m}} \right) - CA_1^{ - 1}B} \right){v_1}\left( \gamma  \right) + {\gamma _{12}}\left( {{I_m} + CA_1^{ - 1}A_2^{ - 1}B} \right){v_2}\left( \gamma  \right) \\&+& {\gamma _{13}}\left( {{I_m} + CA_1^{ - 1}A_3^{ - 1}B} \right){v_3}\left( \gamma  \right) - {\gamma _{12}}{\gamma _{23}}CA_1^{ - 1}A_2^{ - 1}A_3^{ - 1}B{v_3}\left( \gamma  \right) = \alpha _1^*{u_1}\left( \gamma  \right),
\end{eqnarray*}
which can be verified by the first equation of (\ref{setofeq}). Thus, equation (\ref{meq}) holds for the computed vectors $w_i(\gamma), (i=1,2,3)$, concluding that $\lambda_1,\lambda_2$ and $\lambda_3$ are some eigenvalues of $\textnormal{K}_{X_*}$.

\section{$k$ Prescribed eigenvalues}
Keep in mind the matrix K as in (\ref{matrixk}). In this section, improving and generalizing the methodology, necessary definitions and lemmas of previous section and other related works, we aim to extend the results of Section \ref{sec2} to the case of $k\le m$ arbitrary eigenvalues. For the sake of simplicity, we will follow almost the same terminology and writing style, as previous section.
\subsection{Lower bounds for the optimal perturbation}\label{lowek}

For scalars $\gamma  = \left\{ {{\gamma _{12}},{\gamma _{13}}, \ldots {\gamma _{1k}},{\gamma _{23}}, \ldots ,{\gamma _{2k}}, \ldots ,{\gamma _{k - 1k}},} \right\}$, a given set $\Lambda  = \{ \lambda _1 ,\lambda _2 , \ldots ,\lambda _k \}$ and a square matrix $T\in\mathbb{C}^{d\times d}$, define matrix $Q_T(\gamma)$ as 
\[{Q_T}(\gamma ) = \left[ {\begin{array}{*{20}{c}}
  {T - {\lambda _1}{I_l}}&{{\gamma _{12}}{I_l}}& \cdots &{{\gamma _{1k}}{I_l}} \\ 
  0&{T - {\lambda _2}{I_l}}& \ddots & \vdots  \\ 
   \vdots & \ddots & \ddots &{{\gamma _{k - 1k}}{I_l}} \\ 
  0& \ldots &0&{T - {\lambda _k}{I_l}} 
\end{array}} \right],\qquad\gamma  \in \mathbb{C},\]
thus, if $T \in {\mathbb{C}^{l \times l}}$ has $\lambda _1 ,\lambda _2 , \ldots ,\lambda _k$ as some of its eigenvalues, then for all $\gamma  \in \mathbb{C}$ it holds that ${s_{kl-(k-1)} }\left( {Q_T(\gamma )} \right)=0.$ 
\begin{corollary}\label{colorank}
Let $\lambda _1 ,\lambda _2 , \ldots ,\lambda _k$ be some 
eigenvalues of $T \in {\mathbb{C}^{l \times l}}$. Then for all $\gamma  \in \mathbb{C}$  we have $\textnormal{rank}\left( {{Q_T}(\gamma )} \right) \leqslant kl - k.$
\end{corollary}
Let
\[\begin{gathered}
\mathcal{A} = \left[ {\begin{array}{*{20}{c}}
{A - {\lambda _1}{I_n}}&{\gamma_{12} {I_n}}& \ldots &{\gamma_{1k} {I_n}} \\ 
0&{A - {\lambda _2}{I_n}}& \ddots & \vdots  \\ 
\vdots & \ddots & \ddots &{\gamma_{k-1k} {I_n}} \\ 
0& \ldots &0&{A - {\lambda _k}{I_n}} 
\end{array}} \right],{\mathcal{B} = \textnormal{diag}\left[ {\underbrace {B,B, \ldots ,B}_{\textnormal{m times}}} \right],}
\end{gathered} \]
and matrices $\mathcal{X}$ and $\mathcal{C}$ be defined similar to $\mathcal{A}$ and $\mathcal{B}$, respectively. If spectrum of the matrix $ \textnormal {K}_X$ includes the set $\Lambda$, then Corollary \ref{colorank} concludes
\begin{eqnarray*}
\textnormal {rank}\left( \left[{\begin{array}{*{20}{c}}
\mathcal A&\mathcal B \\ 
\mathcal	C&\mathcal X 
\end{array}}\right] \right) &=&\textnormal {rank}\left( {\left[ {\begin{array}{*{20}{c}}
{{K_X } - {\lambda _1}I_{m+n}}&{\gamma_{12} I_{m+n}}& \ldots &{\gamma_{1k} I_{m+n}} \\ 
0&{{K_X } - {\lambda _2}I_{m+n}}& \ddots & \vdots  \\ 
\vdots & \ddots & \ddots &{\gamma_{k-1k} I_{m+n}} \\ 
0& \ldots &0&{{K_X } - {\lambda _m}I_{m+n}} 
\end{array}} \right]} \right) \\&\leqslant& k\left( {n + m - 1} \right),
\end{eqnarray*}
then, by applying Theorem \ref{mainrank3}, $\textnormal {rank}\left( \left[{\begin{array}{*{20}{c}}
\mathcal A&\mathcal B \\ 
\mathcal	C&\mathcal X 
\end{array}}\right] \right) =  \rho \left( \gamma  \right) + \textnormal {rank}\left( {\mathcal{S}_k\left( X ,\gamma\right)} \right),$ in which
\begin{equation}\label{mains}
\begin{gathered}
\mathcal{M}\left( \gamma  \right) = \left( {{I_{nk}} - \mathcal{A}{\mathcal{A}^\dag }} \right)\mathcal{B},\qquad \mathcal{N}\left( \gamma  \right) = \mathcal{C}\left( {{I_{nk}} - {\mathcal{A}^\dag }\mathcal{A}} \right), \hfill \\
\mathcal{S}_k\left( X ,\gamma\right) = \left( {I - \mathcal{N}{\mathcal{N}^\dag }} \right)\left( {\mathcal{X} - \mathcal{C}{\mathcal{A}^\dag }\mathcal{B}} \right)\left( {I - {\mathcal{M}^\dag }\mathcal{M}} \right), \hfill \\ 
\end{gathered}
\end{equation}
\begin{equation*}\label{rho}
\rho\left(\gamma\right)  = \textnormal{rank}\left( {\left[ {\mathcal{A},\mathcal{B}} \right]} \right) + \textnormal{rank}\left( {\left[ \begin{gathered}
\mathcal{A} \hfill \\
\mathcal{B} \hfill \\ 
\end{gathered}  \right]} \right) - \textnormal{rank}\left( \mathcal{A} \right)=kn,
\end{equation*}	
thus, $\rho \left( \gamma  \right) + \textnormal{rank}\left( {{\mathcal{S}_k}\left( {X,\gamma } \right)} \right) \leqslant k\left( {n + m - 1} \right)$ implies $\textnormal{rank}\left( {{\mathcal{S}_k}\left( {X,\gamma } \right)} \right) \leqslant k\left(m-1\right)$
and consequently for every $\gamma\in\mathbb{R}$, it holds that ${s_{k\left( { m - 1} \right) + 1}}\left( {{\mathcal{S}_k}\left( {X,\gamma } \right)} \right) = 0.$ Next lemma gives the desired lower bound.


\begin{lemma}\label{lemmalowerbound}
Assume that the matrix $\textnormal{K}$ as in (\ref{matrixk}) and the set $\Lambda$ of $k\le m$ scalars are given. If for $X\in \mathbb{C}^{m\times m}$, spectrum of $\textnormal {K}_X=\left[ {\begin{array}{*{20}{c}}
A&B \\ 
C&X 
\end{array}} \right]$ includes $\Lambda$ then
\[{\left\| {X - D} \right\|_2}\ge s_{\kappa} \left( \gamma  \right)= {s_{k\left( {m - 1} \right) + 1}}\left( {{\mathcal{S}_k}\left( {D,\gamma } \right)} \right) ,\qquad \gamma  \in \mathbb{C}.\]
\end{lemma}

\subsection{Construction of the optimal perturbation}\label{sectionconstruction}
The main goals considered in this section is finding an optimal perturbed matrix $X_*$ satisfying ${\left\| {X_* - D} \right\|_2}\ge s_{\kappa} \left( \gamma  \right)$ for some value of $\gamma$. Now, at first step, we need to have an explicit formula for the elements of $\mathcal{A}^\dagger$. To this end, note that
\[{\mathcal{A}^\dag } = {\left[ {{\mathrm{a}_{ij}}} \right]_{k \times k}} = \left\{ {\begin{array}{*{20}{c}}
{{0_n}}&{i > j,} \\ 
{{{\left( {A - {\lambda _i}{I_n}} \right)}^{ - 1}}}&{i = j,} \\ 
{ - {{\left( {A - {\lambda _i}{I_n}} \right)}^{ - 1}}\left( {\sum\limits_{r = i + 1}^k {{\gamma_{ir}\mathrm{a}_{rj}}} } \right)}&{i < j.} 
\end{array}} \right.\]
where each block $\mathrm{a}_{ij}$ is an n-by-n square matrix. So, using (\ref{mains}) we have
\[{\mathcal{S}_k}\left( {D,\gamma } \right) = {\left[ {{\mathrm{s}_{ij}}} \right]_{k \times k}} = \left\{ {\begin{array}{*{20}{c}}
{{0_m}}&{i > j,} \\ 
{\left( {D - {\lambda _i}{I_m}} \right) - C{{\left( {A - {\lambda _i}{I_n}} \right)}^{ - 1}}B}&{i = j,} \\ 
{\gamma_{ij} {I_m} + C{{\left( {A - {\lambda _i}{I_n}} \right)}^{ - 1}}\left( {\sum\limits_{r = i + 1}^k {\gamma_{ir}{\mathrm{a}_{rj}}} } \right)B}&{i < j.} 
\end{array}} \right.\]
in which ${\mathrm{s}_{ij}}, \left(i,j \in \left\{ {1,2, \ldots ,k} \right\}\right)$ are $m\times m$ matrices. Assume that vectors
\begin{equation*}
u(\gamma)=\left[ {\begin{array}{*{20}{c}}
{{u_1}(\gamma )}\\
\vdots \\
{{u_k}(\gamma )}
\end{array}} \right], v(\gamma)=\left[ {\begin{array}{*{20}{c}}
{{v_1}(\gamma )}\\
\vdots \\
{{v_k}(\gamma )}
\end{array}} \right] \in {\mathbb{C}^{km}}~({u_j}(\gamma ),{v_j}(\gamma ) \in {\mathbb{C}^{m\times1}},j = 1, \ldots ,k),
\end{equation*}
is a pair of left and right singular vectors of ${s_\kappa }\left( \gamma \right)$, respectively. Let  $u(\gamma)$ and $v(\gamma)$ be unit vectors, define the $m\times k$ matrices
$
U(\gamma) = [{u_1}(\gamma ), \ldots ,{u_k}(\gamma )]$ and $V(\gamma) = [{v_1}(\gamma ), \ldots ,{v_k}(\gamma )]$, and note that 
\begin{equation}\label{skuv}
{\mathcal{S}_k}\left( {D,\gamma } \right)v(\gamma ) = {s_\kappa }\left(\gamma \right)u(\gamma ),
\end{equation}

Some properties of ${s_\kappa }\left( \gamma \right)$ such as being bounded and the fact that if for some value $\gamma$, we have ${s_\kappa }\left( \gamma \right)=0$, then for all $\gamma$ this equation holds, can be derived analogous to proof of Lemma 26 of \cite{graciavelasco}. More importantly, the following lemma, which can be verified by considering Lemma 3.5 of \cite{lipert} and part (iii) of Lemma 26 of \cite{graciavelasco},  concludes that there exists a finite point $\gamma  \in \mathbb{R}$ where the function ${s_\kappa }\left( \gamma \right)$ attains its maximum value.
\begin{lemma}\label{tend20}
If for $i>j$ it holds that rank $\mathrm{s}^k_{ij}\ge k$, then ${s_\kappa }\left( {Q_A(\gamma )} \right) \to 0$ as $\gamma  \to \infty$.
\end{lemma}

Let $\gamma_*>0$ be a point where the singular value ${s_\kappa }\left(\gamma \right)$ attains its maximum value and ${s_\kappa }\left( \gamma_* \right)=\alpha_k^* >0$ be a simple singular value of ${\mathcal{S}_k}\left( {D,\gamma_* } \right)$.
\begin{lemma}\label{lem3k}
Let $\gamma_*$ and $\alpha_k^*$ be as defined above. Then, matrices $U\left( {{\gamma _*}} \right)$ and $V\left( {{\gamma _*}} \right)$ have full rank. Additionally $U{\left( {{\gamma _*}} \right)^*}U\left( {{\gamma _*}} \right) = V{\left( {{\gamma _*}} \right)^*}V\left( {{\gamma _*}} \right).$
\end{lemma}
Now, similar to (\ref{delta}), define
$\Delta_*  =  - {\alpha_k^*}U({\gamma _*}){V }({\gamma _*})^\dag$. Lemma \ref{lem3k} yields $ V{(\gamma )^\dag } V{(\gamma )}=I_k$ and ${\left\| \Delta_*  \right\|_2} = \alpha_k^*$.

For every $\left\{ {{a_1}, \ldots ,{a_t}} \right\} \subseteq \left\{ {1,\ldots ,k} \right\},\left(3\le t\le k\right)$,  let ${P_{{a_1}, \ldots ,{a_t}}} = CA_{{a_1}}^{ - 1}A_{{a_2}}^{ - 1} \ldots A_{{a_t}}^{ - 1}B$ where ${a_1},{a_2}, \ldots ,{a_t}$ is a strictly  increasing sequence of numbers. For instance, we have ${P_{124}} = CA_1^{ - 1}A_2^{ - 1}A_4^{ - 1}B$. Then, keeping in mind (\ref{MNP}), equation (\ref{skuv}) turns into
\begin{equation}\label{sys}
\left\{ {\begin{array}{*{20}{c}}
  \begin{gathered}
  {M_1}{v_1}\left( \gamma  \right) + {\gamma _{12}}{N_{12}}{v_2}\left( \gamma  \right) + \left( {{\gamma _{13}}{N_{13}} - {\gamma _{12}}{\gamma _{23}}{P_{123}}} \right){v_3}\left( \gamma  \right) +  \ldots  \hfill \\
   + \left( {{\gamma _{1k}}{N_{1k}} - {\gamma _{12}}{\gamma _{2k}}{P_{12k}} -  +  \ldots  + {\gamma _{12}}{\gamma _{23}} \ldots {\gamma _{k - 1k}}{P_{12...k}}} \right){v_k}\left( \gamma  \right) = s{u_1}\left( \gamma  \right), \hfill \\ \vspace{.1cm}
\end{gathered}  \\ 
  {{M_2}{v_2}\left( \gamma  \right) + {\gamma _{23}}{N_{23}}{v_3}\left( \gamma  \right) +  \ldots  + \left( {{\gamma _{2k}}{N_{2k}} +  \ldots  + {\gamma _{23}}{\gamma _{34}} \ldots {\gamma _{k - 1k}}{P_{2...k}}} \right){v_k}\left( \gamma  \right) = s{u_2}\left( \gamma  \right),} \\ 
   \vdots \vspace{.05cm} \\ 
  {{M_{k - 1}}{v_{k - 1}}\left( \gamma  \right) + {\gamma _{k - 1k}}{N_{k - 1k}}{v_k}\left( \gamma  \right) = s{u_{k - 1}}\left( \gamma  \right),}\vspace{.3cm} \\
  {{M_k}{v_k}\left( \gamma  \right) = s{u_k}\left( \gamma  \right),} 
\end{array}} \right.
\end{equation}
In order to show that spectrum of ${{\textnormal{K}_{{X_*}}}}$ includes $\Lambda$, consider the matrices
\[W = {\left[ {\begin{array}{*{20}{c}}
  {{w_k}\left( \gamma  \right)}&{{w_{k - 1}\left( \gamma  \right)}}& \ldots &{{w_2}\left( \gamma  \right)}&{{w_1}\left( \gamma  \right)} \\ 
  {{v_k}\left( \gamma  \right)}&{{v_{k - 1}\left( \gamma  \right)}}& \ldots &{{v_2}\left( \gamma  \right)}&{{v_1}\left( \gamma  \right)} 
\end{array}} \right]_{\left( {m + n} \right) \times k}}\]
\[E = {\left[ {\begin{array}{*{20}{c}}
  {{\lambda _k}}&{ - {\gamma _{k - 1k}}}&{ - {\gamma _{k - 2k}}}& \ldots &{ - {\gamma _{1k}}} \\ 
  0&{{\lambda _{k-1}}}& \ddots & \ddots & \vdots  \\ 
   \vdots &0& \ddots &{ - {\gamma _{23}}}&{ - {\gamma _{13}}} \\ 
  {}& \vdots & \ddots &{{\lambda _{2}}}&{ - {\gamma _{12}}} \\ 
  0&0& \ldots &0&{{\lambda _1}} 
\end{array}} \right]_{k \times k}},\]
it is worth noting that in the matrix $E$, elements of every diagonal upper the principle diagonal are reversed and multiplied by -1, in comparison to $\mathcal{A}$. Extending the idea used in previous section, we aim to find $k$ vectors $w_k\left( \gamma  \right)$ satisfying K$_{{X_*}}W = WE$. Clearly, $\lambda_1,\ldots,\lambda_k$ are eigenvalues of upper triangular matrix $E$. It can be easily verified that if K$_{{X_*}}W = WE$ and $\phi_i$ is an eigenvector of $E$ corresponding to $\lambda_i$, then $\left(\lambda_i, W\phi_i\right)$ is an eigenpair of K$_{{X_*}}$. The first row of matrix multiplication K$_{{X_*}}W = WE$, yields the vectors $w_i\left(\gamma\right), \left(i=1,2,\ldots,k\right)$ as follows
\begin{equation}\label{wks}
\begin{array}{*{20}{c}}
  {{w_k}\left( \gamma  \right) =  - A_k^{ - 1}B{v_k}\left( \gamma  \right),} \\ 
  {{w_{k - 1}}\left( \gamma  \right) =  - A_{k - 1}^{ - 1}B{v_{k - 1}}\left( \gamma  \right) + {\gamma _{k - 1k}}A_{k - 1}^{ - 1}A_k^{ - 1}B{v_k}\left( \gamma  \right),} \\ 
   \vdots  \\ 
  {\begin{array}{*{20}{c}}
  {{w_1}\left( \gamma  \right) =  - A_1^{ - 1}B{v_1}\left( \gamma  \right) + {\gamma _{12}}A_1^{ - 1}A_2^{ - 1}B{v_2}\left( \gamma  \right) +  \ldots } \\ 
  { + {{\left( { - 1} \right)}^k}\left( {{\gamma _{12}}{\gamma _{23}} \ldots {\gamma _{k - 1k}}A_1^{ - 1}A_2^{ - 1} \ldots A_k^{ - 1}B} \right){v_k}\left( \gamma  \right)} 
\end{array}} 
\end{array}
\end{equation}
The second row of K$_{{X_*}}W = WE$ also confirms what has obtained for $w_i\left(\gamma\right), \left(i=1,2,\ldots,k\right)$. By way of illustration, multiplying the second row of K$_{X_*}$ and $W$ by the first column of $W$ and $E$, respectively, obtains
\[C{w_k}\left( \gamma  \right) + {X_*}{v_k}\left( \gamma  \right) = {\lambda _1}{v_k}\left( \gamma  \right),\]
considering the fact that ${X_*}{v_i}\left( \gamma  \right) = D{v_i}\left( \gamma  \right) - \alpha _*^k{u_i}\left( \gamma  \right),\left( {i = 1,2, \ldots k} \right)$ and substituting ${w_k}\left( \gamma  \right)$ from (\ref{wks}) into the above equation conclude
\[ - CA_1^{ - 1}B{v_k}\left( \gamma  \right) + D{v_k}\left( \gamma  \right) - \alpha _*^k{u_k}\left( \gamma  \right) = {\lambda _1}{v_k}\left( \gamma  \right),
\]
or equivalently
\[\left( {D - {\lambda _1}{I_m}} \right) - C{\left( {A - {\lambda _1}{I_n}} \right)^{ - 1}}B{v_k}\left( \gamma  \right) = \alpha _*^k{u_k}\left( \gamma  \right),\]
which is indeed the last equation of (\ref{sys}). Performing similar computations for the rest of columns will obtain other equations of (\ref{sys}).

\begin{remark}
If K is a normal matrix, then $\alpha_k^*$ may be repeated, indeed it holds that ${s_{k\left( {m - 1} \right) + 1}}\left( {{\mathcal{S}_k}\left( {D,\gamma_* } \right)} \right)={s_{k\left( {m - 1} \right)}}\left( {{\mathcal{S}_k}\left( {D,\gamma_* } \right)} \right).$ As can be seen from its geometric interpretation shown in Figure \ref{fig:normal}, the graphs of these two singular value functions intersect at $\gamma_*$. Therefore, not only is $\alpha_k^*$ not simple, but $(\gamma_*, \alpha_k^*)$ is a sharp point of ${s_{k\left( {m - 1} \right) + 1}}\left( {{\mathcal{S}_k}\left( {D,\gamma } \right)} \right)$ and consequently this function is not analytic (see Theorem \ref{analytic}). Thus, we may have  $U{\left( {{\gamma_*}} \right)^*}U\left( {{\gamma_*}} \right) \neq V{\left( {{\gamma_*}} \right)^*}V\left( {{\gamma_*}} \right)$. To address this problem, the idea given in \cite{iknormal} can be easily developed.

\begin{figure}
\centering
\includegraphics[width=0.7\linewidth]{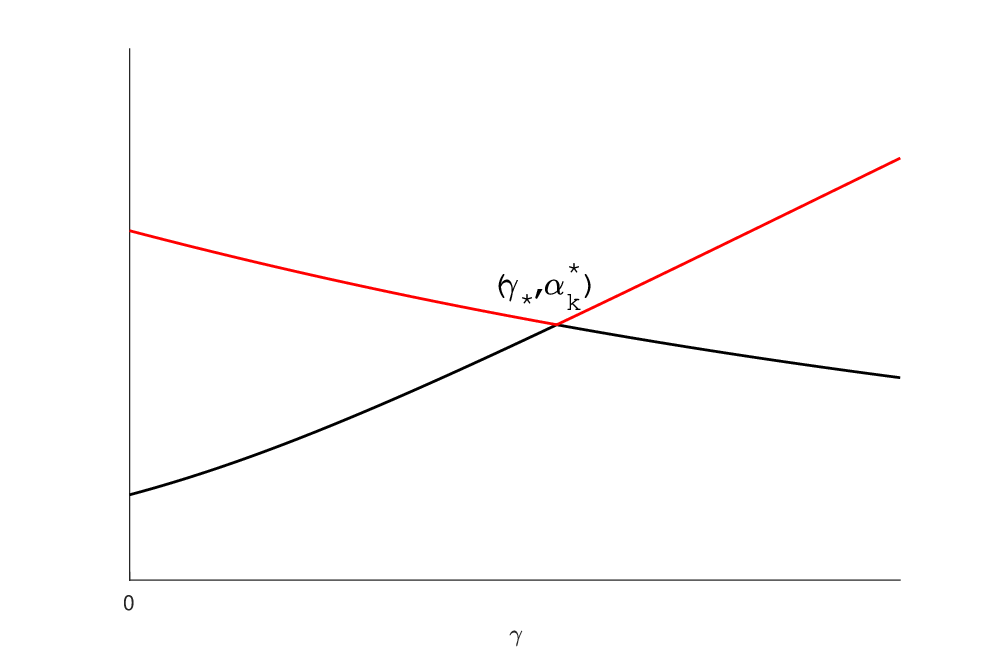}
\caption{Graphs of ${s_{k\left( {m - 1} \right) + 1}}\left( {{\mathcal{S}_k}\left( {D,\gamma_* } \right)} \right)$ (the black graph) and ${s_{k\left( {m - 1} \right)}}\left( {{\mathcal{S}_k}\left( {D,\gamma_* } \right)} \right)$(the red graph) when K is a normal matrix.}
\label{fig:normal}
\end{figure}
\end{remark}

\section{Making optimal change to every main diagonal block to have prescribed eigenvalues}\label{every}
The idea of having desired eigenvalues by making optimal changes to southeast submatrix can be extended to every main diagonal square block. To do this, it is enough to substitute the main diagonal block that should be perturbed for the southeast one, using permutation matrices. To be more precise, assume that a given matrix K$_{d\times d}$ is partitioned such that $A_{11}, A_{22}, \ldots, A_{nn}$ are its main diagonal square blocks (not necessarily of equal size), and $P_{k}$ is the permutation matrix obtained from the $d\times d$ identity matrix by permutation of rows $k$ and $d$.

To construct the optimal perturbation of $A_{kk}$ in order to have prespecified eigenvalues, it is needed to multiply K from the left and right by $P_{k}$. Then, redefine the block matrices $A,B,C$ and $D$ (the submatrix $D$ is now $A_{kk}$ ) to compute ${\mathcal{S}_k}\left( {D,\gamma } \right)$ as described in the beginning of Subsection (\ref{sectionconstruction}). Finally, after making optimal changes to $D(=A_{kk})$, by reusing the permutation matrix $P_k$, return all the blocks in their origin.

For instance, consider $3\times3$ block matrix K of the form 
\begin{equation}\label{matrixk3}
\textnormal{K} = {\left[ {\begin{array}{*{20}{c}}
  {{A_{11}}}&{{A_{12}}}&{{A_{13}}} \\ 
  {{A_{21}}}&{{A_{22}}}&{{A_{23}}} \\ 
  {{A_{31}}}&{{A_{32}}}&{{A_{33}}} 
\end{array}} \right]},\;\;\;\textnormal{in which}\;\;\; A_{11} \in {\mathbb{C}^{n \times n}},{A_{22}} \in {\mathbb{C}^{m \times m}},{A_{33}} \in {\mathbb{C}^{k \times k}}.
\end{equation}

Assume that $A_{22}$ is the only block that can be perturbed. Then, multiplying K from the both sides by ${P_2} = \left[ {\begin{array}{*{20}{c}}
  1&0&0 \\ 
  0&0&1 \\ 
  0&1&0 
\end{array}} \right]$ moves $A_{22}$ to the southeast place. The other blocks move as follows
\[{P_2}\textnormal{K}{P_2} = \left[ {\begin{array}{*{20}{c}}
  {{A_{11}}}&{{A_{13}}}&{{A_{12}}} \\ 
  {{A_{31}}}&{{A_{33}}}&{{A_{32}}} \\ 
  {{A_{21}}}&{{A_{23}}}&{{A_{22}}} 
\end{array}} \right],\]
more importantly, the blocks $A,B,C$ and $D$ are of the form
\[A = \left[ {\begin{array}{*{20}{c}}
  {{A_{11}}}&{{A_{13}}} \\ 
  {{A_{31}}}&{{A_{33}}} 
\end{array}} \right],\;\;B = \left[ {\begin{array}{*{20}{c}}
  {{A_{12}}} \\ 
  {{A_{32}}} 
\end{array}} \right],\;\;C = \left[ {\begin{array}{*{20}{c}}
  {{A_{21}}}&{{A_{23}}} 
\end{array}} \right],\;\;D = {A_{22}},\]
which will be applied to obtain ${\mathcal{S}_k}\left( {D,\gamma } \right)$.
The rest of computations to construct the optimal perturbation $\Delta_*$ is completely similar to what fully explained in Subsections \ref{cons3} as well as \ref{sectionconstruction}.

\section{Conceivable application and numerical examples}
In this section, we review some possible applications of the problem addressed in this paper. Also, the validity and effectiveness of our method and discussions in the previous sections are verified by a comprehensive numerical example. 

\textit{Dynamic systems. }Let us first, borrow some definitions and concepts of \cite{gova,govapaper} about dynamic system that is a system of ordinary differential equations as follows
\begin{equation}\label{ode}
\dfrac{{dx}}{{dt}} = G\left( {x,\alpha } \right),\qquad x \in {\mathbb{R}^n},\alpha  \in {\mathbb{R}^m},G\left( {x,\alpha } \right) \in {\mathbb{R}^n},
\end{equation}
and are ubiquitous in mathematical modeling. See, for instance, \cite{murr,govakhib}. $G$ is function of the state variable $x$ and the parameter $\alpha$. 
The equilibrium solutions of \ref{ode} are those for which $x$ is a constant and $G\left(x,\alpha\right)=0$. The necessary and sufficient condition for structural stability of an equilibrium is that the Jacobian matrix $G_x$ has no eigenvalues on the imaginary axis of the complex plane \cite{govapaper}.
Moreover, a solution to (\ref{ode}) is called a Bogdanov-Takens (BT) point of order $k$ if $G_x$ has exactly an eigenvalue zero with algebraic multiplicity $k$ and geometric multiplicity 1. The role and impact of BT points are well reviewed in section 8.4 of \cite{kuz} and section 7.3 of \cite{guc}.

Now, assume that for a particular given vector $x_0$ and a given dynamic system as in (\ref{ode}), $G_{x_0}$ does not satisfy the condition for structural stability. To cope with the problem and find the nearest matrix to the Jacobian matrix with desired properties,  one can employ the method of this paper to find the closest matrix to $G_{x_0}$ with no eigenvalues of the form $\pm i\omega ,\left( {\omega  > 0} \right)$, then consider it as a reasonable approximation of  $G_{x_0}$. This is also true about a given point to be a BT point.

\textit{Spectral clustering. } The results of this paper can also be applied in spectral clustering. See \cite[Chapter 16]{zaki} and references therein for the theory and applications. To clarify, assume that $G$ is an undirected $d$-regular graph and $V$ and $E$ are  the sets of vertices and edges, respectively. A connected component of $G$ is a maximal set of vertices such that there is an edge between each of them. If $\left|V\right|=n$, then the adjacency matrix of $G$ is defined as follows
\[A = {\left[ {{a_{ij}}} \right]_{n \times n}} = \left\{ {\begin{array}{*{20}{c}}
  1&{\textnormal{if there is an edge between }{v_i}\textnormal{ and }{v_j}} \\ 
  0&{\textnormal{otherwise}} 
\end{array}} \right.\]
considering the diagonal matrix $D = \textnormal{diag}\left\{ {\deg \left( {{v_1}} \right), \ldots ,\deg \left( {{v_n}} \right)} \right\}$ in which deg$(v_i)$ denotes the degree of vertex $v_i$, the Laplacian matrix of $G$ is defined as $L=D-A$.

Spectral theory aims to partition $G$ to $k$ subgraphs with minimum cut (cutting the relation of subgraphs). The related version of this problem becomes:
\[\mathop {\min \textnormal{Trace}\left( {{U^T}LU} \right)}\limits_{{U^T}U = I},\]
The solution of this problem, i.e., $U^*$ is the $k$ eigenvectors of $L$ corresponding to it's $k$ smallest eigenvalues. Also, 
\[\textnormal{Trace}\left( {{U^T}LU} \right) = \sum\limits_{i = n}^{n - k + 1} {{\lambda _i}}, \]
but it is known that for a graph with $k$ distinct subgraphs, we have $\min \textnormal{Trace}\left( {{U^T}LU} \right) = 0,$ which means that ${\lambda _n} = {\lambda _{n - 1}} =  \ldots  = {\lambda _{n - k + 1}} = 0.$ So, $\lambda=0$ has $k$ multiplicity. Therefore, finding $\bar{L}$ such that $\bar{L} \simeq L$ and has zero eigenvalue with multiplicity $k$ is reasonable. This problem can be modeled as
\[\mathop {\min \left\| {\bar L - L} \right\|}\limits_{L \in \Gamma }, \]
in which $\Gamma$ is the set of all matrices of the same size as $L$ that have zero as an eigenvalue of multiplicity$\ge k$.


In particular, eigenvalues of $L$ play a significant role in clustering stability, determining the number of clusters and clustering quality. When it comes to the number of clusters, there is a well known theorem which states that the number of connected components in the graph is equal to the multiplicity of zero as an eigenvalue of $L$. See for example  \cite[Proposition 1.3.7]{spectra}. On the other hand, ideally (neither having noise nor weakly connected components), the number of clusters corresponds to the number of connected components.

Suppose that for a set of given data, having $k$ clusters is desired while 0 is not an eigenvalue of $L$ with algebraic multiplicity $k$. To cope with this problem, we can replace $L$ with its closets approximation that meets the required condition using established method in the paper.  This, would make changes to the current clusters and can be construed as clustering by reverse engineering.

Furthermore, if clustering quality and stability are concerned it can be said that the bigger spectral gap, the more stability of clustering \cite{stability}. The spectral gap is the difference between the two smallest eigenvalues of $L$. Thus, with the aim of improving stability of clustering one can construct an optimal perturbation of $L$ such that the spectral gap of the perturbed matrix is as big as required to guarantee the stability.

It should be noticed that replacing $L$ with an approximation that satisfies the desired conditions, must be performed with great caution ensuring that the modified Laplacian matrix preserves structure of the data both physically and conceptually.

As a numerical experiment, suppose that a matrix K$_{6\times6}$ partitioned as in (\ref{matrixk3}) and a set $\Lambda=\left\{1,2-i,\sqrt{3}\right\}$ are given. Assume that it is desired to find a perturbation of K such that not only just the middle block $A_{22}$ is allowed to change, but also this change is as minimum as possible (with respect to spectral norm), provided that $\Lambda$ is a subset of spectrum of the perturbed matrix. 

Since we are given three prescribed eigenvalues, $A_{22}$ must be a 3-square block, at least. Therefore, considering $n=2, m=3$ and $k=1$, K can be partitioned as
\[
\textnormal{K}=
\left[
\begin{array}{c|c|c}
\begin{array}{*{20}{c}}
  3&{ -5} \\ 
  10&4
\end{array} & \begin{array}{*{20}{c}}
  4&{4}&{-5} \\ 
  {4}&-9&{-6}
\end{array} & \begin{array}{*{20}{c}}
  {4} \\ 
  {7} 
\end{array}\\
\hline
\begin{array}{*{20}{c}}
  {-3}&{-10} \\ 
  6&{2} \\ 
  4&-2 
\end{array} & \begin{array}{*{20}{c}}
  9&{-2}&-4 \\ 
  {-10}&-1&6 \\ 
  -1&6&-1
\end{array} & \begin{array}{*{20}{c}}
  { - 10} \\ 
  -7\\ 
  {5} 
\end{array}\\
\hline
\begin{array}{*{20}{c}}
  -1&{-7} 
\end{array} & \begin{array}{*{20}{c}}
  { -3}&{ 2}&{ - 6} 
\end{array} & 5
\end{array}
\right].
\]

By adopting the procedure described in the paper and using \textit{fminsearch} MATLAB function that is a multidimensional unconstrained nonlinear minimization, it can be seen that the function $s_3(\gamma)$ attains its maximum value at
\[{\gamma _*} = \left\{ {\begin{array}{*{20}{c}}
  {{\gamma _*}_{12} = 5.7459,}&{{\gamma _*}_{23} =  1.2047,}&{{\gamma _*}_{13} = -9.5084 -0.9111i} 
\end{array}} \right\},
\]
and $\alpha _3^* = 4.9119$. It is also worth mentioning that $\alpha_3^*$ is an isolated singular value of ${\mathcal{S}_3}\left( {D,\gamma_* } \right)$, as needed. Additionally, two matrices
\[U = \left[ {\begin{array}{*{20}{c}}
  {{{ - 0}}{{.0462  +  0}}{{.0022}}i}&{{{0}}{{.2930  -  0}}{{.0773}}i}&{{{ - 0}}{{.1985  +  0}}{{.0054}}i} \\ 
  {{{0}}{{.1694  -  0}}{{.0160}}i}&{{{0}}{{.4330  -  0}}{{.0922}}i}&{{{ - 0}}{{.4832  +  0}}{{.0987}}i} \\ 
  {{{0}}{{.2684  -  0}}{{.0456}}i}&{{{0}}{{.5515  -  0}}{{.0892}}i}&{{{ - 0}}{{.1097  +  0}}{{.0114}}i}  
\end{array}} \right],\]
and
\[V = \left[ {\begin{array}{*{20}{c}}
  {{{ - 0}}{{.1307}}}&{{{0}}{{.0804  -  0}}{{.0554i}}}&{{{ - 0}}{{.0481  -  0}}{{.0021i}}} \\ 
  {{{ - 0}}{{.2446  +  0}}{{.0113i}}}&{{{ - 0}}{{.5444  +  0}}{{.0294i}}}&{{{0}}{{.0848  -  0}}{{.0209i}}} \\ 
  {{{0}}{{.1661  -  0}}{{.0236i}}}&{{{0}}{{.5233  -  0}}{{.1391i}}}&{{{ - 0}}{{.5258  +  0}}{{.0918i}}} 
\end{array}} \right],\]
satisfy what claimed in Lemma \ref{lem3}, more precisely, are full rank and
\[{\left\| {U{{\left( {{\gamma _*}} \right)}^*}U\left( {{\gamma _*}} \right) - V{{\left( {{\gamma _*}} \right)}^*}V\left( {{\gamma _*}} \right)} \right\|_2} = {{2}}{{.9065}} \times {{1}}{{{0}}^{ - 6}}.\]

\begin{figure}[t]
\centering
\includegraphics[width=0.6\linewidth]{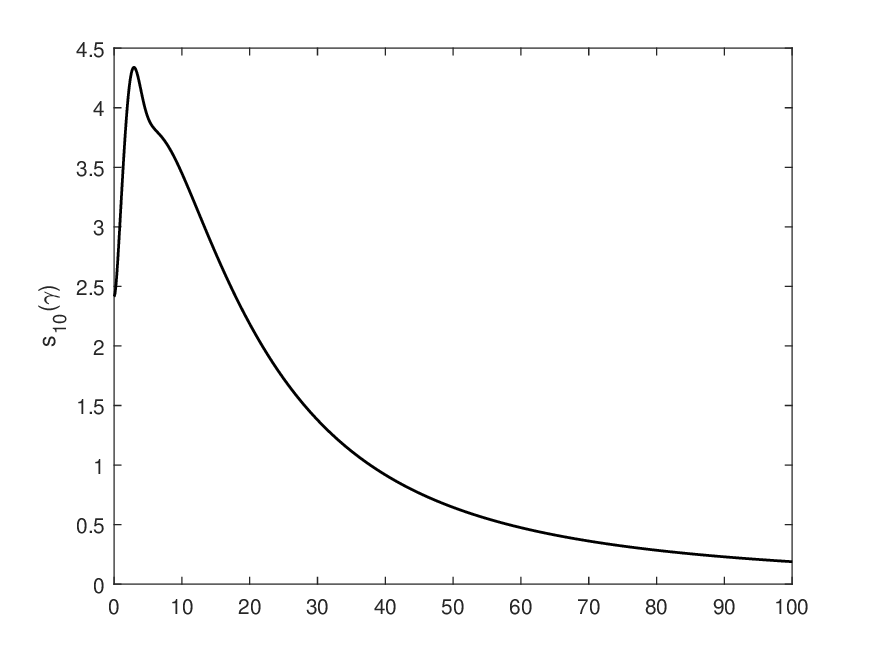}
\caption{The graph of the function $s_{10}\left(\gamma\right)$.}
\label{fig:s10}
\end{figure}

As shown in Figure \ref{fig:s10}, $s_{10}\left(\gamma\right)\to 0$, when $\left| \gamma  \right| \to \infty$, as proved. Furthermore, the optimal perturbation $\Delta_*$ as in (\ref{delta}) is of the form
\[\Delta_*  = \left[ {\begin{array}{*{20}{c}}
  {{{ - 4}}{{.6844  +  0}}{{.1673i}}}&{{{0}}{{.6878  -  0}}{{.0261i}}}&{{{ - 1}}{{.2805  -  0}}{{.2018i}}} \\ 
  {{{1}}{{.2046  -  0}}{{.2466i}}}&{{{ - 0}}{{.4433  +  0}}{{.3773i}}}&{{{ - 4}}{{.7156  +  0}}{{.1950i}}} \\ 
  {{{0}}{{.7989  -  0}}{{.0625i}}}&{{{4}}{{.7596  -  0}}{{.8124i}}}&{{{ - 0}}{{.3161  -  0}}{{.2655i}}} 
\end{array}} \right],\]
and spectrum of the perturbed matrix is as follows
\[\sigma \left( {{\textnormal{K}_{{X^*}}}} \right) = \left\{ {\begin{array}{*{20}{c}}
  { - 0.7516 + 12.0308i,}&{ - 0.5027 - 11.5211i,}&{10.0784 + 0.7695i,} \\ 
  {1,}&{2 - i,}&{1.7321\left( { = \sqrt 3 } \right)} 
\end{array}} \right\}.\]

\end{document}